\documentclass[twoside,bezier]{amsart}
\usepackage[left=2cm,right=2cm]{geometry}

\usepackage{amsmath,amsthm,amsfonts,amssymb,amscd,euscript}
\usepackage{setspace}

\usepackage[figuresright]{rotating}
\newtheorem{thm}{Theorem}[section]
\newtheorem{cor}[thm]{Corollary}
\newtheorem{lem}[thm]{Lemma}

\theoremstyle{definition}\newtheorem{defn}[thm]{Definition}
\theoremstyle{remark}
\newtheorem{rem}[thm]{Remark}
\numberwithin{equation}{section}
\newcommand{\A}{\mathcal{A}}

\newcommand{\filebegin}{\begin{document}}
\newcommand{\fileend}{\end{document}}
\makeatother
\def\thefootnote{}
\newcommand{\lo}{\longrightarrow}
\newcommand{\NMM}{\hspace*{2mm}}

\usepackage[left=2cm,right=2cm]{geometry}

\def\n{\noindent}%

\numberwithin{equation}{section}
\def\mapdown#1{\Big\downarrow\rlap
{$\vcenter{\hbox{$\scriptstyle#1$}}$}}
\begin{document}

\vspace*{2cm}
\begin{center}
{\bf\large On an equivalence of the topology of algebraic cone metric spaces
and metric spaces}
 \\[0.5cm]
{Saeedeh Shamsi Gamchi, Mohammad Janfada, Assadollah Niknam\\
Department of Pure Mathematics, Ferdowsi University of Mashhad, Mashhad, P.O. Box 1159-91775, Iran\\

Email: saeedeh.shamsi@gmail.com\\
Email: mjanfada@gmail.com\\
Email: dassamankin@yahoo.co.uk} \\[2mm]

\end{center}%
\vspace*{0.5cm}
\begin{quotation}
\noindent
{\footnotesize
{\sc Abstract.}
In this paper, we prove that the topology induced by
algebraic cone metric coincides with the topology induced by the
metric obtained via a nonlinear scalarization function, i.e. any algebraic cone metric space is metrizable. Furthermore, the notion of algebraic cone normed space is introduced and also normability of the topology of this space is proved.}
\end{quotation}

\n \textbf{2010 Mathematics subject classification: } Primary: 06A11, 46A19; Secondary 47H10.%\\

\markboth
{S. Shamsi Gamchi , M. Janfada, A. Niknam}
 {On an equivalence of algebraic cone metric spaces
and metric spaces}

%%%%%%%%%%%%%%%%%%%%%%%%%%%%%%%%%%%%%%%%%%%%%%%%%%%%%%%%%%%%%%%%%%%%%%%%%%%%%%%%%%%%%%%%%%%%%%%%%%%%%%%%%%%%%%%%%%

%%%%%%%%%%%%%%%%%%%%%%%%%%%%%%%%%%%%%%%%%%%%%%%%%%%%%%%%%%%%%%%%%%%%%%%%%%

\section {\sc Introduction}
The concept of metric and any concept related to metric play a very important role not only in pure mathematics but
also in other branches of science involving mathematics especially in computer science, information science, and biological
science.\\Replacing the real
numbers, as the codomain of metrics, by ordered Banach spaces
we obtain a generalization of metric spaces. Such generalized
spaces called cone metric spaces, were introduced by Rzepecki
\cite{R}. Long-Guang and Xian \cite{HZ} announced the notion of a cone metric
space by replacing real numbers with an ordering Banach space which is the same as the definition of Rzepecki. After their work many authors
attempt to adjust the theory of cone metric to ordinary metric
space, by proving the most important standard results on fixed
point theory and functional analysis. Recently, Erdal Karapinar \cite{KAR}
studied fixed point theorems in cone Banach spaces, and Abdeljawad
et al. \cite{DA} studied some properties of cone Banach spaces. In
\cite{SC}, Sonmez and Cakalli studied the main properties of cone
normed spaces and
proved some results in cone normed spaces and complete cone normed spaces.\\Ordered normed spaces and cones have applications in applied
mathematics and optimization theory \cite{DN}. A useful approach for analyzing a vector optimization problem is to reduce
it to a scalar optimization problem. Nonlinear scalarization functions play an
important role in this reduction in the context of non-convex vector optimization
problems. Recently this has been applied by Du \cite{WSD} to investigate
the equivalence of vectorial versions of fixed point theorems of contractive mappings in
TVS-cone metric spaces and scalar versions of fixed point theorems in general
metric spaces in usual sense. The authors in \cite{cak} proved that the topology induced by
topological vector space valued cone metric coincides with the topology induced by the
metric obtained via a nonlinear scalarization function. Those results based on nonlinear scalarization function obtained by using a topology in $E$. In this paper we remove this condition by study only algebraic case. Our aim is to generalize the notion of nonlinear scalarization function to prove that any algebraic cone metric space with its inherited topology is metrizable, i.e. the topology induced by the algebraic cone metric coincides with the topology induced by an appropriate metric. Also normability of the algebraic cone normed spaces with their induced topologies is discussed. Note that in \cite{nshj}, we introduced the notion of algebraic cone metric and it has been shown that every algebraic cone metric space has a Hausdorff topology.
%========================================================
%===================================================

\section {\sc Preliminaries}

In this section we recall the concept of algebraic cone metric
space and some of its elementary properties which have been studied in \cite{nshj}.
\begin{defn}
\cite{nk} Let $E$ be a real vector space and
$P$ be a convex subset of $E$. A point $x\in P$ is said to be an
algebraic interior point of $P$ if for each $v\in E$ there exists
$\epsilon>0$ such that $x+tv\in P$, for all $t\in [0,\epsilon]$.
\end{defn}
Note that the above definition is equivalent to the following statement:\\A point $x$ is an algebraic interior point of the convex set $P\subset E$ if $x\in P$ and for each $v\in E$ there exists $\epsilon>0$ such that $[x, x+\epsilon v]\subset P$, where $[x, x+\epsilon
v]=\{\lambda x+ (1-\lambda)(x+\epsilon v): \forall \lambda \in [0,1]\}$.

The set of all algebraic interior points of $P$ is called its
algebraic interior and is denoted by $aint P$. Moreover, $P$ is
called algebraically open if $P=aint P$. Equivalently, $P$ is
called algebraically open if its intersection with every
straight line in $E$ is an open interval (possibly empty). For
example every convex open set in $\mathbb{R}^d$ is algebraically
open. \\

Suppose that $E$ is a real vector space  with its
zero vector $\theta$ and $P\subset E$ is a
non-empty set such that $P+P\subset P$, $\lambda P\subset P$ $(\lambda\geq0)$, $P\cap(-P)=\{\theta\}$. In this case we will say that $P$ is an algebraic
cone in $E$. For a given algebraic cone $P$ in $E$, a
partial ordering $\preceq _{a}$ on $E$ with respect to $P$ is defined
by $x\preceq _{a} y$ if and only if $y-x\in P$. Furthermore, we write $x\ll_{a} y$ whenever $y-x\in aintP$ and  we say that $(E,P)$ is an algebraic cone space. If for each $x\in E$ and $y\in P\setminus \{\theta\}$ there exists $n\in \mathbb{N}$ such that $x\preceq_{a} ny$, we say that $(E,P)$ has the Archimedean property.\\
One can easily see that $P=\{(x_{1},x_{2},...,x_{n})\in\mathbb{R}^{n}: x_{i}\geq 0, ~ i=1,2,..., n)\}$ is an algebraic cone in $\mathbb{R}^{n}$ such that $(\mathbb{R}^{n},P)$ has the Archimedean property. Also $P=\{f\in C_{\mathbb{R}}[a,b]: f(x)> 0, ~ \forall x\in [a,b]\}\cup \{\theta\}$ is an algebraic cone which $(C_{\mathbb{R}}[a,b],P)$ has the Archimedean property. But there exists a real vector space with an algebraic cone which does not have the Archimedean property. For example, in the real vector space $C_{\mathbb{R}}(0,\infty)$ if we consider the algebraic cone $P=\{f\in C_{\mathbb{R}}(0,\infty): f(x)\geq 0, ~ x\in (0,\infty)\}$ then it is easy to see that $P$ does not have the Archimedean property. Indeed, if $f(x)=x^{2}$ and $g(x)=\frac{1}{x}$ then there is no $n\in\mathbb{N}$ such that $g(x)\leq nf(x)$, for all $x\in (0,\infty)$.\\

\begin{lem}\label{c}
Let $(E,P)$ be an algebraic cone space and $aint P\neq\emptyset$. Then\\(i) $P+aintP\subset aintP$.\\(ii) $\alpha aintP\subset aintP$, for each scalar $\alpha >0$.\\(iii) For any $x,y,z\in X$, $x\preceq_{a} y $ and $y\ll_{a} z$ imply that $x\ll_{a}z$.
\begin{proof} Let $x\in aintP, ~ y\in P $ and $v$ be an arbitrary element in $E$. By definition of algebraic interior point, there exists $\epsilon>0$ such that $x+tv\in P$, for each $t\in [0, \epsilon ]$. We have $(x+y)+tv\in P$, for each $t\in [0, \epsilon ]$, since $P+P\subset P$ and $y\in P$. This prove (i).\\
(ii) Let $x\in aintP$, $\alpha >0$ and $v$ be an arbitrary element in $E$. By definition of algebraic interior point, there exists $\epsilon>0$ such that $x+t\frac{v}{\alpha}\in P$, for all $t\in [0,\epsilon]$, hence $\alpha x+tv \in P$, for all $t\in [0,\epsilon]$.\\
(iii) is trivial by using (i).
\end{proof}
\end{lem}
\begin{defn}
Let $(E,P)$ be an algebraic cone space, $aint P\neq\emptyset$ and $d_{a}:X\times X\rightarrow E$ be a
vector-valued function that satisfies:\\(ACM1) For all $x,y\in X$, such that $x\neq y$,
$\theta\ll _{a} d_{a} (x,y)$ and $d_{a}(x,y)=\theta$ if and
only if $x=y$,\\(ACM2) $d_{a}(x,y)=d_{a}(y,x)$ for all $x,y\in
X$,\\(ACM3) $d_{a}(x,y)\preceq _{a} d_{a}(x,z)+d_{a}(z,y)$  for
all $x,y,z\in X$.\\Then $d_{a}$ is called an algebraic cone metric on $X$ and $(X,d_{a})$ is said to be an algebraic cone metric space.
\end{defn}
\begin{thm}\label{n}
Let $(X,d_{a})$ be an algebraic cone metric
space. Then the collection $\{B_{a}(x,c) : c\in aintP,~ x\in X \}$ forms a subbasis for a Hausdorff topology on $X$, where $$B_{a}(x,c):=\{y\in X:d_{a}(x,y)\ll_{a} c\}.$$
\begin{proof} Trivially $\bigcup_{x\in X, c\in aintP} \ B_{a}(x,c)=X$, so the collection $\{B_{a}(x,c) : c\in aintP,~ x\in X \}$ forms a subbasis for a topology on $X$. Now we show that this topology is Hausdorff. Let $x,y\in X$ and $x\neq y$, take $\theta\ll_a c=d_{a}(x,y)$. The facts that $P\cap (-P)=\{\theta\}$ and $d_{a}$ has the property (ACM3) imply that $B_{a}(x,\frac{c}{3})\cap B_{a}(y,\frac{c}{3})=\emptyset$. Therefore the topology induced by the above collection is Hausdorff.
\end{proof}
\end{thm}

\section {\sc Main results }

In this section we suppose that $aintP\neq\emptyset$ and $(E,P)$ be an algebraic cone space. We prove the metrizability of the induced topology on an algebraic cone metric spaces. But before proceeding further, we need a couple of useful lemmas:
\begin{lem}\label{l}
Let $e\in aintP$ . Then $$E=\{\lambda e-c ~: ~ \lambda\in \mathbb{R}^{+}\setminus \{0\}, ~ c\in aint P\}.$$
\begin{proof} Let $e\in aint P$ and $z$ be an arbitrary element of $E$. By definition of algebraic interior point for $v=-z$ there exists $\epsilon>0$ such that $e+tv\in P$, for each $t\in [0,\epsilon]$. Thus for $t =\frac{\epsilon}{2}$, we have $e+\frac{\epsilon}{2}v\in P$ and hence $\frac{2}{\epsilon}e-z\in P$ i.e. $z \preceq_{a} \frac{2}{\epsilon} e$. By Lemma \ref{c}(ii),(iii), we may assume that $z \ll_{a} \lambda e$, for some $\lambda\in \mathbb{R}^{+}\setminus \{0\}$. Thus the proof is complete.
\end{proof}
\end{lem}

\begin{lem}\label{a}
Let $e\in aintP$, then we have\\(i) If $z\in \lambda e-P$ for some $\lambda \in \mathbb{R}$, then for each $\mu > \lambda$, $z\in \mu e-aintP$, in particular, for each $\mu > \lambda$, $z\in \mu e-P$.\\(ii) For any $z\in E$, there exists a real number $\alpha\in \mathbb{R}$ such that $z\notin \alpha e-P$.
\begin{proof}
(i) Let $\mu > \lambda$ and $z\in \lambda e-P$. We have$$\mu e-z=(\mu - \lambda)e+\lambda e -z\subset aintP+P\subset aintP.$$Thus, $z\in \mu e-aintP\subset \mu e-P$.\\(ii) Let us assume that there exists $z_0\in E$ such that $z_0\in P_{\lambda ,e}$, for all $\lambda\in \mathbb{R}$. Then from (i), $z_0\in \lambda e-aintP$, for all $~\lambda\in\mathbb{R}$. Thus,$$\{\lambda e-z_0 :~\lambda\in\mathbb{R} \}\subset aintP$$equivalently,$$\{-\lambda e-z_0 :~\lambda\in\mathbb{R} \}\subset aintP.$$Now by Lemma \ref{l}, for each $y\in E $, there exist $c\in aintP$ and $\alpha\in \mathbb{R}^{+}\setminus \{0\}$ such that $-y=\alpha e-c$. Hence,
\begin{align*}
y&=-\alpha e+c\\
&=(-\alpha e-z_0)+z_0+c\\
&\in aintP+z_0+aintP\subset z_0+aintP.
\end{align*}
Thus,$$E\subset z_0+aintP.$$This contradicts $aintP\neq E$.
\end{proof}
\end{lem}
The nonlinear scalarization function $\xi_{e}: E\rightarrow
\mathbb{R}$ is
defined as follows:$$\xi_{e}(y)=\inf \{r\in\mathbb{R}:  y\in re-P\}$$for all $y\in E$, where $E$ is a topological vector space and $P$ is a closed convex cone in $E$ such that $P \cap(-P)=\{\theta\}$ (see \cite{chy}). The original version is due to Gerstewitz \cite{Gerst}. Its first appearance in English seems to be due to Luc \cite{Luc}.\\
In this approach real vector spaces are used as the domain of the nonlinear scalarization function, instead of topological vector spaces.\\For any $y\in E$, put $$M_{e,y}=\{r\in \mathbb{R}: y\in re-P\}. $$By applying Lemma \ref{l}, it is easy to see that $M_{e,y}$ is non-empty. Moreover, it is bounded below. Indeed, we can assume that, for each $r\in\mathbb{R}$, there exists $\lambda_{r}\in\mathbb{R}$ such that $\lambda_{r}< r$ and $y\in \lambda_{r}e-P$. By Lemma \ref{a}(ii), there exists $\alpha\in \mathbb{R}$ such that $y\notin \alpha e-P$.
Hence, Lemma \ref{a}(i) implies that $y\notin \mu e-P$, for each $\mu <\alpha$ which is a contradiction. Thus $M_{e,y}$ is bounded below.
\begin{defn}
Let $E$ be a real vector space and $P$ be an algebraic cone in $E$. For a given $e\in aint P$, the nonlinear scalarization function is defined by:$$\xi_{e}(y)=\inf M_{e,y}.$$
\end{defn}
\begin{lem}\label{b}
For any $e\in aint P$, the function $\xi_{e}$ has the following properties:\\
1) $\xi_{e}(\theta) = 0$.\\
2) $\xi_{e}(e)=1$.\\
3) $y\in P$ implies $\xi _{e}(y)\geq0$.\\
4) $\xi_{e}(y) <r$ if and only if $ y\in re-aintP$.\\
5) if $y_1 \preceq_{a} y_{2}$, then $\xi_{e}(y_1)\leq \xi_{e}(y_{2})$, for each $y_1 ,y_2 \in E$.\\
6) $\xi_{e}$ is subadditive on $E$.\\
7) $\xi_{e}$ is positively homogeneous on $E$ (i.e. $\xi_{e}(\lambda y)=\lambda \xi_{e}(y)$, for each $y\in E$).\\
8) $\xi_{e}(y)>0$ for each $y\in aintP$.

\begin{proof}For any $y\in E$ we have$$M_{e,y}=\{r\in \mathbb{R}: y\preceq_{a} re\}.$$
1) Since $e\in P$ and $P$ is a cone in $E$, so $re\in P$, for each $r\geq 0$. On the other hand $P\cap (-P)=\{\theta\}$ thus  $M_{e,\theta}=[0,\infty)$.\\
2) Let $r\in M_{e,e}$. Then we have $e\preceq_{a}re$, i.e. $(r-1)e\in P$. Since $P\cap(-P)=\{\theta\}$, so $r\geq 1$. Thus $\xi_{e}(e)=1$.\\
3) Let $y\in P$ and $r\in M_{e,y}$. Then we have $\theta \preceq_{a}y\preceq_{a} re$, hence $\theta \preceq_{a} re$, i.e. $r\in M_{e,\theta}$. Thus $r\geq 0$.\\
4) Let $\xi_{e}(y)<r$. Then there exists $r'<r$ such that $y\preceq_{a}r'e$, on the other hand by Lemma\ref{c}(ii), we have $r'e\ll_{a}re$. Hence, Lemma\ref{c}(iii) implies that $y\ll_{a} re$. Conversely, let $y\in re-aintP$. By definition of algebraic interior point, for $v=-e$, there exists $\epsilon>0$ such that $(re-y)+tv\in P$, for all $t\in [0,\epsilon]$. Hence, there exists $r'<r$ such that $r'e-y\in P$. Thus $\xi_{e}(y)\leq r'<r$.\\
5) Let $y_1 ,y_2 \in E$ and $y_1 \preceq_{a} y_{2}$. It is easy to see that $M_{e,y_{2}}\subset M_{e,y_{1}}$. Thus $\xi_{e}(y_1)\leq \xi_{e}(y_{2})$.\\
6) Let $y_1 ,y_2 \in E$, $r_1 \in M_{e,y_{1}}$ and $r_2 \in M_{e,y_{2}}$. Then $y_1 +y_2 \preceq_{a}(r_1 +r_2)e$, hence $\xi_{e}(y_1 +y_2)\leq r_1 +r_2$. By taking infimum on $r_1\in M_{e,y_{1}} $ and then on $r_2\in M_{e,y_{2}}$ we have $$\xi_{e}(y_1 +y_2)\leq \xi_{e}(y_1 )+\xi_{e}(y_2).$$
7) Follows from definition of $\xi_{e}$.\\
8) Since $y\in aintP$, so by definition of algebraic interior point for $v=-e$, there exists $\epsilon>0$ such that $y+tv\in P$, for all $t\in [0,\epsilon]$, hence $y-\frac{\epsilon}{2} e \in P$, i.e. $\frac{\epsilon}{2} e \preceq_{a} y$. Therefore, properties (2), (5) and (7) imply that $\xi_{e}(y)>0$.
\end{proof}
\end{lem}
By using the idea of Wie-shie Du \cite{WSD}, we can assert the following theorem.
\begin{thm}\label{f}
Let $(X,d_{a})$ be an algebraic cone metric space and $e\in aintP$. Then $d_e: X\times X \rightarrow [0,\infty)$ defined by $d_e=\xi_{e}\circ d_a$ is a metric.
\begin{proof}
By (ACM1), (ACM2) and Lemma \ref{b}, $d_e(x,y)\geq 0$ and $d_e(x,y)=d_e(y,x)$ for all $x, y \in X$. If $x =y$, then by (ACM1),
$d_e(x,y)=\xi_{e}(\theta)=0$. Conversely, if $d_e(x,y)=0$, then from Lemma \ref{b}(8) and (ACM1) we have $d_{a}(x,y)=\theta$, consequently $x=y$. Applying (5) and (6) of Lemma \ref{b} and (ACM3), $d_e$ satisfies the triangle inequality which implies that $d_e$ is a metric.
\end{proof}
\end{thm}
\begin{lem}\label{e}
Let $(X, d_{a})$ be an algebraic cone metric space, $e\in aintP$, $d_e=\xi_{e}\circ d_a$, and $x\in X$. Then $B_{d_{e}} (x, r) = B_{a}(x,re)$, where $B_{d_{e}} (x,r)=\{y\in X: d_e(x,y)<r\}$ and $B_{a}(x,re)=\{y\in X:d_{a}(x,y)\ll_{a} re\}$.
\begin{proof}
Let $y\in B_{d_{e}}(x,r)$. Then $\xi_{e}(d_a (x,y))<r$, so by Lemma \ref{b}(4), we have $d_a(x,y)\in re-aintP$ i.e. $d_a(x,y)\ll_a re$, hence $y\in B_a(x,re)$. Moreover, by Lemma \ref{b}(4), it is easy to verify that $B_a(x,re)\subset B_{d_{e}} (x,r)$. Thus the proof is complete.
\end{proof}
\end{lem}

\begin{rem}\label{m}
Let $(X,d_{a})$ be an algebraic cone metric space. It is proved in Theorem \ref{n} that the collection $\{B_{a}(x,c) : c\in aintP,~ x\in X \}$ forms a subbasis for a Hausdorff topology on $X$, where $B_{a}(x,c):=\{y\in X:d_{a}(x,y)\ll_{a} c\}$. Consider the metric defined by $d_e=\xi_{e}\circ d_{a}$. Let $\tau _{d_e}$ and $\tau _{a}$ denote the topology induced by the metric $d_e$, and the
topology induced by the algebraic cone metric $d_{a}$, respectively. By using Lemma \ref{e}, one can see $\tau _{d_e}\subset \tau _{a}$.
\end{rem}
In the sequel we assume that $(E,P)$ has the Archimedean property.\\Let $(X, d_{a})$ be an algebraic cone metric space. Then, for each $c_1,c_2 \in aintP$, there exists $c\in aintP$ such that
$c \ll_{a}c_1$ and $c \ll_{a} c_2$. Indeed, by Lemma \ref{c}(ii),(iii) and the Archimedean property of $E$, there exists $n\in\mathbb{N}$ such that $\frac{c_1}{n}\ll_{a} c_2$ and $\frac{c_1}{n}\ll_{a} c_1$. Take $c=\frac{c_1}{n}$, so we find $c \ll_{a}c_1$ and $c \ll_{a} c_2$.\\

Note that this implies that the collection $\{B_{a}(x,c) : c\in aintP,~ x\in X \}$ forms a basis for $\tau_{a}$ the topology of X which is induced by $d_{a}$. Consequently, if $x\in X$ and $\{x_{n}\}$ is a sequence in $X$, then $\{x_{n}\}$ converges to $x$ with respect to $\tau_{a}$ if and only if for every $\theta\ll_{a} c$ there exists a natural number $N$ such that for all $n>N$, $d_{a}(x_{n}, x)\ll_{a} c$.
\begin{defn}
Let $(X, d_{a})$ be an algebraic cone metric space and $\{x_{n}\}$ be a sequence in $X$. Then\\(i) $\{x_{n}\}$ is a Cauchy sequence whenever for every $\theta\ll_{a} c$ there exists a natural number $N$ such that for all $m,n>N$, $d_{a}(x_{n}, x_{m})\ll_{a} c$.\\(ii) $(X, d_{a})$ is said to be a complete algebraic cone metric space if every Cauchy sequence is convergent.
\end{defn}

Let $(X,d_a)$ be an algebraic cone metric space, $x\in X$ and $\{x_{n}\}$ be a sequence in X. Let $d_e$ be the same as in Theorem \ref{f}, then it is easy to see that\\
(i) $\{x_n\}$  converges to x if and only if $d_e(x_n,x)\rightarrow 0$ as $n\rightarrow \infty$;\\
(ii) $\{x_n\}$ is a Cauchy sequence in $(X,d_a)$ if and only if $\{x_n\}$ is a Cauchy sequence (in usual sense) in $(X,d_e)$;\\
(iii) $(X,d_a)$ is complete algebraic cone metric space if and only if $(X,d_e)$ is a complete metric space.

\begin{cor}\label{h}
Let $(X, d_{a})$ be a complete algebraic cone metric space. Suppose that a mapping $T :X\rightarrow X$ satisfies the contractive condition$$d_{a}(T x,Ty) \preceq_{a}\lambda d_{a}(x, y),$$for all $x,y\in X$, where $\lambda\in (0, 1)$ is a constant. Then T has a unique fixed point in X. Moreover, for each $x\in X$, the iterative sequence ${T^{n} x}$ converges to the fixed point.
\begin{proof}
Put $d_e=\xi_{e}\circ d_a$. Applying the assertion befor this corollary, $(X,d_e)$ is a complete metric space. By Lemma \ref{b}(4), (7), we have$$d_{a}(T x,Ty) \preceq_{a}\lambda d_{a}(x, y) ~ \Rightarrow ~ d_e(T x,Ty) \leq\lambda d_e(x, y).$$
Therefore, the conclusion follows from the Banach contraction principle (see \cite{skc}).
\end{proof}
\end{cor}
Following a similar argument as in Corollary \ref{h}, one can easily obtain vectorial versions of Kannan's fixed point theorem,
Chatterjea's fixed point theorem and others (see \cite{reza,rhoad,abbas}). Moreover, the equivalence between scalar versions and vectorial
versions of these results can also be established easily.

\begin{thm}\label{k}
Let $(X, d_a)$ be an algebraic cone metric space. Then there exists a metric on $X$ which induces the same topology on $X$ as
the topology induced by $d_a$.
\begin{proof}
Let $e\in aintP$ and consider the metric defined by $d_e=\xi_{e}\circ d_{a}$. Let $\tau _{d_e}$ and $\tau _{a}$ denote the topology induced by the metric $d_e$, and the
topology induced by the algebraic cone metric $d_{a}$, respectively. By Remark \ref{m}, we have $\tau_{d_e}\subset\tau_{a}$. It is sufficient to show that $\tau_{a}\subset\tau_{d}$. Now let $U\in \tau_{a}$, then there exists $c\in aintP$ such that $B_a(x,c)\subset U$. Since $E$ has the Archimedean property, so there exists $n\in \mathbb{N}$ such that $e\ll_{a} nc$. Hence $B_a(x,\frac{e}{n})\subset B_a(x,c)\subset U$. On the other hand, by Lemma \ref{e}, $B_d(x,\frac{1}{n})= B_a(x,\frac{e}{n})$. Therefore, $B_d(x,\frac{1}{n})\subset U$ which implies that $U\in\tau _{d}$. Thus the topology $\tau_{d}$ coincides with the topology $\tau_{a}$.
\end{proof}
\end{thm}
Note that this theorem shows that there are many metrics on $X$ whose topologies are the same as the topology induced by $d_a$.\\
Now, we give a definition of algebraic cone normed space.
\begin{defn}\label{i}
Let $X$ be a vector space over $F$ ($\mathbb{R}$ or $\mathbb{C}$) and $\|.\|_{a}:X\rightarrow E$ be a
mapping that satisfies:\\(ACN1)  $\theta\ll_a \|x\|_{a}$ for all $x\in
X\setminus\{\theta_X\}$ and $\|x\|_{a}=\theta$ if and only if $x=\theta_{X}$, where
$\theta_{X}$ is the zero vector in $X$,\\(ACN2)  $\|\alpha x\|_{a}=|\alpha|\|x\|_{a}$
for all $x\in X$ and $\alpha \in F$,\\(ACN3)  $\|x+y\|_{a}\preceq_a \|x\|_{a}+\|y\|_{a}$.\\Then $\|.\|_{a}$ is called an algebraic cone norm on $X$ and $(X,\|.\|_{a})$ is called an algebraic cone normed space.
\end{defn}
\begin{lem}\label{j}
Let $(X,\|.\|_{a})$ be an algebraic cone normed space. Then for any $e\in aintP$, $\|.\|_{e}: X \rightarrow [0,\infty)$ defined by $\|.\|_e:=\xi_e\circ \|.\|_a$ is a norm on $X$.
\begin{proof}
It follows by Lemma \ref{b} and Definition \ref{i}.
\end{proof}
\end{lem}

\begin{thm}
Every algebraic cone normed space is normable, i.e. the topology of the algebraic cone normed space coincides with the topology of
the norm defined in Lemma \ref{j}.
\begin{proof}
It is clear that $d_a(x,y)=\|x-y\|_a$ is an algebraic cone on $X$. Put $\|.\|_e:=\xi_e\circ \|.\|_a$, and $d_e(x,y)=\|x-y\|_e= \xi_e \circ\|x-y\|_a$, hence we have $d_e(x,y)=\xi_e\circ d_a(x,y)$. Then it follows from Lemma \ref{j} and Theorem \ref{k} that the topology
of the algebraic cone normed space $(X, \| .\|_a)$ coincides with the topology of the normed space $(X,\|.\|_e)$.
\end{proof}
\end{thm}
\section {\sc Conclusion}
We would like to emphasize that the present work contains not only the proofs of metrizability of an algebraic cone metric space and normability of an algebraic cone normed space, but also some useful theorems cause to generalize nonlinear scalarization functions to cosider some optimization problems.

%{\bf Acknowledgments.} The authors would like to sincerely thank the referee for very
%valuable comments and suggestions.

\providecommand{\bysame}{\leavevmode\hbox
to3em{\hrulefill}\thinspace}

%%%%%%%%%%%%%%%%%%%%%%%%%%%%%%%%%%%%%%%%%%%%%%%%%%%%%%%%%%%%%%%%%%%%%%%%%

\end{document}